\documentclass[11pt, reqno]{amsart}

\usepackage{fullpage}
\usepackage{amsmath}
\usepackage{amsfonts}
\usepackage{amsthm}
\usepackage{amssymb}
\usepackage{enumitem}
\usepackage{hyperref}
\usepackage{microtype}

\usepackage[backgroundcolor=lightgray]{todonotes}

\theoremstyle{plain}
\newtheorem{theorem}{Theorem}

\newtheorem{proposition}[theorem]{Proposition}

\theoremstyle{definition}

\newtheorem{example}[theorem]{Example}

\theoremstyle{remark}
\newtheorem{remark}[theorem]{Remark}

\newcommand{\beqn}{\begin{eqnarray*}}
	\newcommand{\eeqn}{\end{eqnarray*}}

\usepackage{mathtools}

\DeclarePairedDelimiterX{\dset}[2]{\lbrace}{\rbrace}{#1\;\delimsize|\;#2}

\newcommand{\N}{\mathbb{N}}
\newcommand{\R}{\mathbb{R}}

\newcommand{\ball}{{\overline{B}}}

\title[Yet another approach to measuring distance in hyperspaces]{Yet another approach to measuring distance in some hyperspaces}

\author[D.Bugajewski]{Dariusz Bugajewski}
\address{D. Bugajewski, Department of Nonlinear Analysis and Applied Topology\\
  Faculty of Mathematics and Computer Science\\
  Adam Mickiewicz University, Pozna\'n\\
  ul.\ Umultowska 87\\
  61-614 Pozna\'n\\
  Poland}
\email{ddbb@amu.edu.pl}

\author[P.Ma\'ckowiak]{Piotr Ma\'ckowiak} 
\address{P. Ma\'ckowiak, Department of Nonlinear Analysis and Applied Topology\\
  Faculty of Mathematics and Computer Science\\
  Adam Mickiewicz University, Pozna\'n\\
  ul.\ Umultowska 87\\
  61-614 Pozna\'n\\
  Poland}
\email{piotr.mackowiak@amu.edu.pl}

\date{\today}

\begin{document}

\begin{abstract} 
In this article, we define a new partial distance on hyperspaces consisting of closed balls, which is a partial metric in the sense of Matthews. We prove that these partial metric spaces are complete if and only if the metric spaces from which the corresponding ball spaces arise are complete. Finally, we compare the topology induced by the new partial distance with that induced by the classical Hausdorff distance. 
\end{abstract}

\maketitle
	
\section{Introduction}

Given a metric space $X$ denote by $H(X)$ the family of all nonempty closed and bounded subsets of $X$. Probably the most known metric defined on the family $H(X)$ is the so-called Hausdorff metric (more information about the Hausdorff metric the interested reader can found e.g. in the monographs \cite{AB} and \cite{N}). A lot of articles has been devoted to the study of hyperspaces endowed with the Hausdorff distance. That fact is connected, in particular, with their numerous applications e.g. in digital topology (see e.g. \cite{B}, \cite{BIK}, \cite{HKR} and \cite{S}). Let us mention here a recently published article \cite{HZ} in which it was proven that under $CH$, the Hausdorff metric topology on the hyperspace $CL(X)$ of all nonempty closed subsets of a metric space $(X, d)$ is completely metrizable if and only if $(X, d)$ is completely metrizable. 

A particular place among hyperspaces play hyperspaces consisting of balls included in a given metric space. There various approaches to measuring distance in these hyperspaces. Besides the classical Hausdorff distance one can use some smetrics as e.g. in \cite{ST} and \cite{NT} (let us add also that in the paper \cite{NT} one can find some applications of ball spaces in fixed point theory) or ultrametrics as it was done in the papers \cite{BKK1} and \cite{BKK2}. 
	
Let us now turn our attention to partial metric spaces. Although these spaces are generally non-Hausdorff, many metric-like tools and concepts can be naturally extended to this setting. In particular, as shown in~\cite{BMW}, in partial metric spaces compactness and sequential compactness are equivalent. Initially, partial metric spaces have found applications in denotational semantics of programming languages. However, they also find interesting applications in other areas, such as the geometry of normed spaces (see~\cite{ROSP2}), the domain of words, and complexity spaces (see~\cite{ROSP}). Finally, one has to emphasize  that the fixed point theory in partial metric spaces is quite well developed (see, for instance,~\cite{BM, BMW, HRS, IPR} and the references therein). 

The main goal of this paper is define a new partial distance on hyperspaces consisting of closed balls, which is a partial metric in the sense of Matthews. We  compare the topology induced by the new partial distance with that induced by the classical Hausdorff distance and we examine some basic topological spaces of that partial metric spaces.

The paper is organized as follows. In Section 2, we recall some basic definitions and facts related to partial metric spaces and weighted smetric spaces that will be used in this paper. In Section 3, we introduce a new function $p$ on hyperspaces consisting of closed balls and prove that it is a partial metric in the sense of Matthews. We show that the topologies generated by the Hausdorff metric and that new partial metric are not comparable. We also establish a necessary and  sufficient condition for the completeness of that hyperspace endowed with the partial metric $p$ and we show that an analogue of the Blaschke Theorem on the compactness in the Hausdorff metric does not hold in that case, too.

\section{Preliminaries} 
In this section, we collect some basic notions and results, especially those related to partial metric spaces and weighted semi-pseudo metrics. Further information and examples regarding partial metric spaces can be found, for instance, in~\cite{BW} while more details concerning weighted partial metric spaces can found, for example, in \cite{ST} 

Let $U$ be a non-empty set. A mapping $p \colon U\times U \to [0,+\infty)$ is called a \emph{partial metric} (or \emph{pmetric}) on $U$, if for all $x,y,z\in U$ the following four conditions are satisfied:
\begin{enumerate}[label=\textup{(P\arabic*)}]
 \item $x=y$ if and only if $p(x,x)=p(y,y)=p(x,y)$,
 
 \item $p(x,x)\leq p(x,y)$,

 \item $p(x,y)=p(y,x)$,

 \item $p(x,y)\leq p(x,z)+p(z,y)-p(z,z)$.
\end{enumerate}
The pair $(U,p)$ is called a \emph{partial metric space}, and $p(x,x)$ is referred to as the \emph{size} of $x$. Clearly, any metric space is a partial metric space with all points of the size $0$.

Recall that a sequence $(x_n)_{n\in \mathbb{N}}$ of elements of a partial metric space $(U,p)$ is said to \emph{converge to an element $x\in U$} if
\[
\lim_{n\rightarrow \infty}p(x_n,x)=p(x,x);
\]
it is said to \emph{properly converge} to $x\in U$, if it converges to $x$ and $\lim_{n\to \infty}p(x_n,x_n)=p(x,x)$. Further, a sequence $(x_n)_{n\in \mathbb{N}}$ of elements of a partial metric space $(U,p)$ is said to be a \emph{Cauchy sequence} if there exists (and is finite)
\[
\lim_{n,m\rightarrow\infty}p(x_n,x_m).
\]
Finally, a partial metric space $(U,p)$ is said to be \emph{complete} if every Cauchy sequence $(x_n)_{n\in \mathbb{N}}$ in $(U,p)$ of elements of $U$ is properly convergent.

To every partial metric space $(U,p)$ one can naturally assign metric spaces, which preserve some properties of $(U,p)$. Namely, let for all $x,y\in U$, define metrics $p^m$ and $d_m$ by
\begin{align*}
  p^m(x,y)&:=2p(x,y)-p(x,x)-p(y,y),\\
	d_m(x,y)&:=\max\bigl\{p(x,y)-p(x,x),\ p(x,y)-p(y,y)\bigr\}.
\end{align*}
For example, it is well-known that the completeness of a partial metric space $(U,p)$ is equivalent to the completeness of both metric spaces $(U, p^m)$ and $(U, d_m)$ (see, for example,~\cite[p.~194]{M} or \cite{C}, Prop. 5.24). 

In a partial metric space $(U,p)$, the set $\dset{y \in U}{p(x,y)\leq r}$ may be empty for some $r>0$ and $x \in U$. To avoid this issue, we define the closed ball centered at $x\in U$ with radius $r>0$ as
\[
\ball_p(x,r):=\dset{y \in U}{p(x,y)\leq p(x,x)+r}.
\]
When $p$ is a metric, this definition clearly agrees with the usual notion of a closed ball. 

Let us pass on to semi-pseudo metrics. A mapping $d \colon U\times U \to [0,+\infty)$ is called a \emph{semi-pseudo metric} (or \emph{smetric}) on $U$, if for all $x,y\in U$ the following two conditions are satisfied:
\begin{enumerate}\label{def:smetric}
 \item[\rm{(S1)}] $d(x,x)=0$,
 \item[\rm{(S2)}] $d(x,y)=d(y,x)$.
\end{enumerate}
A weighted smetric space is a triple $(U,d,w)$, where $d$ is an smetric on $U$ and $w:U\to [0,+\infty)$ is a function such that for all $x,y,z\in U$ the following condition is satisfied: 
\begin{enumerate}\label{def:smetric}
 \item[\rm{(S3)}] $d(x,z)\leq d(x,y) + w(y) + d(y,z)$.
\end{enumerate} 
Values of the function $w$ are said to be weights of suitable points of the space $U$.\\
Given a metric space $(X,\rho)$, by $\ball(x,\epsilon)$ we shall denote the closed ball centered at $x\in X$ with radius $\epsilon\geq 0$ (we simply write $\ball$ if we do not indicate the center and the radius of that ball) and, by $B(X)$ - the set of closed balls of $X$, endowed with the smetric $d_s$ defined as follows: 
\[
d_s(\ball_1,\ball_2):=\inf\{\rho(x,y)\ | \ x\in \ball_1,y\in \ball_2\}, \quad \mbox{where} \quad \ball_1, \ball_2 \in B(X). 
\] 

It was proven in \cite{ST} (Proposition 5.2 (1)) that if $X$ is totally convex, then, taking as the weight of any ball its diameter, $B(X)$ becomes a weighted smetric space. 

Observe, that we allow the radius of a ball $\ball\in B(X)$ to be zero; thus, all singletons of $X$ are members of $B(X)$. 

Denote by $B_{d_s}(\overline{A},\epsilon)$ an open ball in $B(X)$ with center at $\overline{A}\in B(X)$ and radius $\epsilon>0$, that is, 
\[
B_{d_s}(\overline{A},\epsilon):=\{\ball \in B(X)\colon d_s(\overline{A},\ball)<\epsilon \}.
\]
In general, the family of such balls is not a basis for any topology, as shows, in particular, the following 

\begin{example}\label{ex:d_s}
Let us consider $(\mathbb R,\rho)$, where $\rho$ denotes the Euclidean metric. Let $A_1=\ball (-1,\frac{1}{2}), A_2=\ball (1,\frac{1}{2}), A_3=\ball (0,\frac{1}{2})\in B(X)$. It is clear that for every $\epsilon>0$, it holds: 
\[
A_3\in B_{d_s}(\overline{A}_1,\epsilon)\cap B_{d_s}(\overline{A}_2,\epsilon).
\] 
Fix $\epsilon=\frac{1}{2}$. For every $\epsilon'\in (0,\frac{1}{4})$, we have 
\[
B_{d_s}(\overline{A}_3,\epsilon')\not\subset B_{d_s}(\overline{A}_2,\frac{1}{2})\cap B_{d_s}(\overline{A}_3,\frac{1}{2}). 
\] 
Indeed, we have 
\[
\ball (-1,\frac{1}{2}+\epsilon')\in B_{d_s}(\overline{A}_3,\epsilon'), 
\]
however 
\[
d_s(\ball (-1,\frac{1}{2}+\epsilon'),A_2)=\inf\left\{|x-y|\colon x\in \left[-\frac32-\epsilon',-\frac{1}{2}+\epsilon'\right], y\in\left[\frac{1}{2},\frac{3}{2}\right]\right\}>\frac{1}{2}, 
\] 
so  
\[
\ball (-1,\frac{1}{2}+\epsilon')\not\in B_{d_s}(\overline{A}_2,\frac{1}{2}). 
\] 
Thus, the family of open balls in smetric space $(B(X),d_s)$ is not a basis for any topology. 
\end{example}

\section{A partial metric defined on ball-sets}

Let $(X,\rho)$ be a metric space. For the sake of simplicity of presentation, we assume that in the considered metric space $(X,\rho)$ any closed ball - as a set - is uniquely identified through the center and the radius of the ball, that is, $\ball(x,r)=\ball(x',r')$ if and only if $x=x'$ and $r=r'$. 

Let $p: B(X)\times B(X)\to [0,+\infty)$ be a mapping defined as follows:  

\begin{equation}\label{pmetric}
 	p(\ball_1,\ball_2):=\left\{
	\begin{array}{lr}
		2r_1,& \text{ for } \ball_1=\ball_2,\\
	d_s(\ball_1,\ball_2)+2r_1+2r_2,& \text{ otherwise},
\end{array}
\right.
\end{equation}
where a $\ball_i\in B(X)$ has center $x_i$ and radius $r_i$ (we use this convention in what follows). 

The crucial property of the above defined function $p$ is given in the following 

\begin{proposition}\label{prop:pm} 
The function $p$ defined by formula (\ref{pmetric}) is a partial metric on $B(X)$. 
\end{proposition} 

\begin{proof} 
The conditions $\rm{(P1)-(P3)}$ are satisfied in an obvious way. Let us check the condition $\rm{(P4)}$. Consider only the case when $\ball_1,\ball_2,\ball_3\in B(X)$ are pairwise distinct. Then for any $v,v'\in \ball_3(z,r_3)$, we have
\begin{align*}
p(\ball_1,\ball_2)&=d_s(\ball_1,\ball_2)+2r_1+2r_2=\inf\{\rho(x,y)\ | \ x\in \ball_1,y\in \ball_2\}+2r_1+2r_2\\
&\leq\rho(x,y)+2r_1+2r_2\leq \rho(x,v)+\rho(v,v')+\rho(v',y)+2r_1+2r_2\\
&\leq\rho(x,v)+2r_1+2r_3+\rho(v',y)+2r_2+2r_3-2r_3.
\end{align*} 
Taking suitable infima on both sides of the above inequality we get
\[
p(\ball_1,\ball_2)\leq p(\ball_1,\ball_3)+p(\ball_3,\ball_2)-p(\ball_3,\ball_3). 
\]
\end{proof} 

\begin{remark} 
Let us notice that the topologies generated by the metric $d_H$ and the partial metric $p$ defined at the beginning of this section need not be comparable. 

Indeed, let us consider again $(\mathbb R,\rho)$, where $\rho$ denotes the Euclidean metric. Denote by $\tau_H$ the topology generated by the Hausdorff metric and, by $\tau_p$ the topology generated by the partial metric $p$ defined above. Moreover, denote by $B_H(\overline{A},\epsilon)$ an open ball in $(B(X),d_H)$ with center at $\overline{A}\in B(X)$ an radius $\epsilon>0$ and analogously, by $B_p(\overline{A},\epsilon)$ in the space $(B(X),p)$. Let $\overline{A}=\ball (0,1)\in B(X)$. Obviously $\overline{A}\in B_{d_H}(\overline{A},\frac{1}{4})$ and $p(\overline{A},\overline{A})=2$. Let $\overline{A'}\in B_p(\overline{A},\epsilon)$, $\overline{A}\neq \overline{A'}$ and $\epsilon>0$. Thus  
\[
\inf\{|x-y|\colon \ x\in \overline{A},y\in \overline{A'}\}+2r_{\overline A}+2r_{\overline{A'}}<\epsilon +2r_{\overline A}, 
\] 
where $r_{\overline A}$ and $r_{\overline{A'}}$ denote the radii of the balls $\overline{A}$ and $\overline{A'}$, respectively. Then $r_{\overline{A'}}<\frac{\epsilon}{2}$. Now, if $\epsilon<\frac{1}{4}$ and $\overline{A'}\in B_p(\overline{A},\epsilon')$, then $\epsilon'\geq \frac{3}{4}$, which means that $\overline{A}\not\in B_{d_H}(\overline{A},\frac{1}{4})$. Hence, the topology $\tau_p$ is not stronger than the topology $\tau_H$. 

Now, fix $\epsilon \in (0,\frac{1}{4})$. Obviously $\overline{A}\in B_p(\overline{A},\frac{1}{4})$. Moreover, if $\overline{A'}\in B_{d_H}(\overline{A},\epsilon)$, $\overline{A'}\neq \overline{A}$, then $|r_{\overline{A}}-r_{\overline{A'}}|\leq \epsilon$, so $r_{\overline{A'}}\in (\frac{3}{4},\frac{5}{4})$, because $\overline{A}=1$. Therefore 
\[
p(\overline{A},\overline{A'})\geq 2+\frac{3}{2}=p(A,A)+\frac{3}{2}, 
\] 
so $\overline{A'}\not\in B_p(\overline{A},\frac{1}{4})$, which implies that $\tau_H$ is not stronger than the topology $\tau_p$. 
\end{remark}

To the partial metric $p$ from Proposition \ref{prop:pm} assign the metric $p^m$ as described in Section 2. Namely, let for $\ball_1, \ball_2\in B(X),$
$$	p^m(\ball_1,\ball_2):=\left\{
	\begin{array}{lr}
		0,& \text{ for } \ball_1=\ball_2,\\
	2d_s(\ball_1,\ball_2)+2r_1+2r_2,& \text{ otherwise},
\end{array}
\right.$$
\begin{remark}
It is clear that the partial metric space $(B(X),p)$ is complete if and only if the metric space $B(X),p^m$ is complete. Moreover, let us notice the topology generated by the metric $p^m$ is stronger than the topology generated by the classical Hausdorff metric $d_H$. It is a consequence of the simple inequality: 
\[
d_H(\ball_1,\ball_2)\leq p^m(\ball_1,\ball_2) \quad \mbox{for any} \quad \ball_1,\ball_2\in B(X).
\] 
\end{remark} 

The following theorem gives a necessary and sufficient condition for the completeness of the space $(B(X),p)$.
\begin{theorem}\label{th:compl}
The partial metric space $(B(X),p)$ is complete if and only if the metric space $(X,\rho)$ is complete.
\end{theorem} 

\begin{proof}
First, assume that the partial metric space $(B(X),p)$ is complete and let 
\[
B^0:=B^0(X):=\{\ball(x,0)\colon x\in X\}=\{\{x\}: x\in X\}\subset B(X). 
\] 
Obviously the mapping $id\colon (B^0,p)\to (X,\rho)$ is an isometry, that is $p(\{x\},\{y\})=\rho(x,y)$ for every $x,y\in X$ and $B^0$ is a closed subspace of the complete metric space $B(X),p^m)$. Hence $(X,\rho)$ is complete. 

To prove the inverse implication assume now that $(X,\rho)$ is complete. Let $(\ball_{n})_{n\in\mathbb N}$ be a Cauchy sequence of balls in $(B(X),p^m)$. It follows immediately that $\lim_{n\to \infty}r_n=0$ and $\lim_{n,m\to \infty}\rho(x_n,x_m)=0$, provided the sequence $(\ball_n)_{n\in\mathbb N}$ is not eventually constant. If the sequence $(\ball_n)_{n\in\mathbb N}$ is eventually constant, then there is nothing to prove. Therefore assume that the sequence $(\ball_n)_{n\in\mathbb N}$ is not eventually constant. By the completeness of $X$, there exists $x_0\in X$ such that $\lim_{n\to \infty}x_n=x_0$. Let $\ball_0:=\ball (x_0,0)$. We have 
\[
p^m(\ball_n,B_0)\leq 2\rho(x_n,x_0)+2r_n\to 0 \quad \mbox{as} \quad n\to \infty. 
\]
The space $(B(X),p^m)$ is complete, so the space $(B(X),p)$ is complete, too. 
\end{proof} 

Now, let us consider the following 

\begin{example}
Let $X:=\R$ be endowed with the Euclidean metric. For every $n\in \N$ define $\ball_n:=\ball(1/n, 1/2)\in B(X)$. We shall show that the sequence $(\ball_n)_{n\in \N}$ has no convergent subsequence in $(B(X),p)$, where $p$ is defined by formula (\ref{pmetric}).
Suppose that $\ball:=\ball(x,r)\in B(X)$, where $x\in X$, $r\geq 0$, is a point in $B(X)$ which is a limit point of the sequence $(\ball_n)_{n\in \N}$. Notice that $p(\ball,\ball)=2r$ and, for all $n\in \N$, except at most one, $p(\ball_n,\ball)\geq\inf_{x\in \ball_n, y\in \ball}|x-y|+2r+1/2$. Thus for a subsequence of $(\ball_n)_{n\in \N}$, converging to $\ball$, (if we slightly abuse notation) we would have $\lim_{n\to \infty}p(\ball_n,\ball)=p(\ball,\ball)$ and, consequently, $2r\geq\lim_{n\to\infty}(\inf_{x\in \ball_n, y\in \ball}|x-y|+2r+1/2)\geq 2r+1/2$, which is impossible. Therefore, there is no convergent subsequence of the sequence $(\ball_n)_{n\in \N}$.
\end{example}
Observe that, as the sequence of balls $(\ball_n)_{n\in \N}$ given in the example above is bounded in the partial metric space $(B(X),p)$, an analogue of the Blaschke Theorem on the compactness in the Hausdorff metric does not hold for the partial metric defined in \eqref{pmetric}.

\end{document}